\documentclass[runningheads]{llncs}
\usepackage{graphicx}
\usepackage{booktabs}
\usepackage{array}
\usepackage{multirow}
\usepackage{algorithm}
\usepackage{algpseudocode}
\usepackage{amssymb}
\usepackage{wrapfig}

\usepackage{amsmath,amsfonts,bm}

\def\eqref#1{equation~\ref{#1}}

\def\1{\bm{1}}

\def\vb{{\bm{b}}}
\def\vc{{\bm{c}}}
\def\vd{{\bm{d}}}

\def\vg{{\bm{g}}}
\def\vh{{\bm{h}}}

\def\vt{{\bm{t}}}

\def\vx{{\bm{x}}}
\def\vy{{\bm{y}}}
\def\vz{{\bm{z}}}

\def\evg{{g}}
\def\evh{{h}}

\def\evx{{x}}
\def\evy{{y}}

\def\mT{{\bm{T}}}

\def\mW{{\bm{W}}}

\DeclareMathAlphabet{\mathsfit}{\encodingdefault}{\sfdefault}{m}{sl}
\SetMathAlphabet{\mathsfit}{bold}{\encodingdefault}{\sfdefault}{bx}{n}

\def\sL{{\mathbb{L}}}

\def\sN{{\mathbb{N}}}

\def\sR{{\mathbb{R}}}
\def\sS{{\mathbb{S}}}

\def\sU{{\mathbb{U}}}
\def\sV{{\mathbb{V}}}

\def\sX{{\mathbb{X}}}

\usepackage{algorithm}
\usepackage{algpseudocode}

\begin{document}
\title{Optimization Over Trained Neural Networks: Taking a Relaxing Walk}
%
%
\author{Jiatai Tong\inst{1} \and
Junyang Cai\inst{1,2} \and
Thiago Serra\inst{1}}
\authorrunning{J. Tong et al.}
%
\institute{Bucknell University, Lewisburg PA, United States 
\email{\{jt037,jc092,thiago.serra\}@bucknell.edu}
\and
University of Southern California, Los Angeles CA, United States
\email{caijunya@usc.edu}}

\maketitle              
\begin{abstract}
Besides training, mathematical optimization is also used in deep learning 
to model and solve formulations over trained neural networks for purposes such as 
verification, compression, and optimization with learned constraints.
However, solving these formulations soon becomes difficult as the network size grows 
due to the weak linear relaxation and dense constraint matrix.  
We have seen improvements in recent years with cutting plane algorithms, reformulations, and an heuristic based on Mixed-Integer Linear Programming (MILP).
In this work, 
we propose a more scalable heuristic based on exploring global and local linear relaxations of the neural network model. 
Our heuristic is competitive with a state-of-the-art MILP solver and the prior heuristic while producing better solutions with increases in input, depth, and number of neurons.

\end{abstract}
\section{Introduction}

There is a natural role for mathematical optimization in machine learning with training, 
where discrete optimization has a ``discreet'' but growing presence in classification trees~\cite{shioda2007crio,dunn2017oct,verwer2017trees,verwer2019trees,hu2019sparse,verhaeghe2020trees,zhu2020multivariate,gunluk2021trees,carrizosa2021survey,aghaei2021trees,demirovic2022murtree,alston2022trees}, decision diagrams~\cite{hu2022diagrams,florio2023diagrams}, decision rules~\cite{dash2018rules,lawless2023rules}, and neural networks~\cite{icarte2019cp,kurtz2021binarized,rosenhahn2022hopfield,patil2022training,bienstock2023principled,thorbjarnarson2023training,bernardelli2023bemi,aspman2023training}.

Now a new role has emerged with predictions from machine learning models being used as part of the formulation of optimization problems. 
For example, imagine that we train a neural network on historical data for approximating an objective function that we are not able to represent explicitly. 
Overall, we start with one or more trained machine learning models, 
and then we formulate an optimization model which---among other things---represents the relationship between decision variables for the inputs and outputs of those trained models. 
Since other discrete decision variables and constraints may be part of such optimization models, 
gradient descent is not as convenient here as it is for training. 

These formulations are \emph{neural surrogate models} if involving neural networks. 
Neural networks are essentially nonlinear, and thus challenging to~model~in~mathematical~optimization. 
However, 
we can use Mixed-Integer Linear Programming~(MILP) for popular activations  
such as the Rectified Linear Unit~(ReLU) \cite{hahnloser2000origin,nair2010rectified,glorot2011rectifier,lecun2015nature,ramachandran2018pop}. 
Neural networks with ReLUs represent piecewise linear functions~\cite{pascanu2013on,arora2018understanding}, which we model in MILP with binary decision variables 
altering the slopes~\cite{serra2018bounding}. 
As with other activations~\cite{cybenko1989approximation,funahashi1989approximate,hornik1989approximator}, 
ReLU networks have been shown to be universal function approximators 
with one hidden layer but enough neurons~\cite{yarotsky2017relu} and with limited neurons per layer but enough layers~\cite{lu2017expressive,hanin2017approximating,park2021width}. 

Many frameworks to formulate neural surrogate models have emerged---JANOS~\cite{bergman2022janos}, OMLT~\cite{ceccon2022omlt}, OCL~\cite{fajemisin2023ocl}, OptiCL~\cite{maragno2023mixed}, and Gurobi Machine Learning~\cite{gurobi2023ml}---in addition to stochastic and robust optimization variants~\cite{dumouchelle2022neur2sp,dumouchelle2023neur2ro,kronqvist2023alternating2sp}. 
The applications in machine learning include 
network verification~\cite{cheng2017resilience,anderson2019strong,anderson2020strong,rossig2021verification,strong2021global}, 
network pruning~\cite{serra2020lossless,serra2021compression,elaraby2023oamip},  
counterfactual explanation~\cite{kanamori2021counterfactual}, 
and constrained reinforcement learning~\cite{delarue2020rlvrp,burtea2023safe}. 
In the broader line of work often denoted as \emph{constraint learning}, 
these models have been used for scholarship allocation~\cite{bergman2022janos}, patient survival in chemotherapy~\cite{maragno2023mixed}, power generation~\cite{murzakhanov2022powerflow} and voltage regulation~\cite{chen2020voltage} in power grids, boiling point optimization in molecular design~\cite{mcdonald2023molecular}, and automated control of industrial operations in general~\cite{say2017planning,wu2020scalable,yang2021control}. 

However, these models can be difficult to solve as they grow in size.  
They have weak linear relaxations due to  
the dense constraint matrix within each layer 
and the big M constraints for each neuron, 
which sparked immediate and continued interest in calibrating big M coefficients~\cite{cheng2017resilience,fischetti2018constraints,liu2021algorithms,tsay2021partition,badilla2023tradeoff} as well as 
in strengthening the formulation and generating cutting planes~\cite{anderson2019strong,anderson2020strong,tsay2021partition}. 
Other improvements include 
identifying stable neurons~\cite{tjeng2019evaluating,xiao2019training}, 
exploiting the dependency among neural activations~\cite{botoeva2020efficient,serra2020empirical}, 
and inducing sparser formulations by network pruning~\cite{say2017planning,cacciola2023pruning}. 
But at the rate of one binary variable per ReLU, 
typically-sized neural networks entail considerably large MILPs, 
hence limiting the applications where these models are solvable within reasonable time.

We may expect that improving scalability will require algorithms exploiting the model structure. For example, 
Fischetti and Jo~\cite{fischetti2018constraints} first observed that a feasible solution for the MILP mapping from inputs to outputs of a single neural network is immediate once a given input is chosen. 
This strategy has been shown effective at least twice~\cite{serra2021compression,perakis2022optimizing}, 
whereas finding a feasible solution for an MILP is generally NP-complete~\cite{conforti2014integer}.
Another example of special structure comes from ReLU networks representing piecewise linear functions. 
Within each part of the domain mapped as a linear function, which is denoted as a \emph{linear region}, there is a direction for locally improving the output. 
In fact, Perakis and Tsiourvas~\cite{perakis2022optimizing} developed a local search heuristic that moves along adjacent linear regions by solving restrictions of the MILP model with some binary variables fixed. 
However, the reliance on MILP eventually brings scalability issues back---although much later in comparison to solving the model without restrictions.

But is there hope for optimization over linear regions at scale? 
That is akin to thinking about MILPs as unions of polyhedra in disjunctive programming~\cite{balas2018dp}. 
While earlier studies have shown that the number of linear regions may grow fast on model dimensions~\cite{pascanu2013on,montufar2014on,telgarsky2015benefits,raghu2017expressive,arora2018understanding}, 
later studies have shown that there are architectureal tradeoffs limiting such growth~\cite{montufar2017notes,serra2018bounding,serra2020empirical,cai2023pruning}. 
Moreover, the networks with typical distributions of parameters have considerably fewer linear regions~\cite{hanin2019complexity,hanin2019deep}; and  
gradients change little between adjacent linear regions~\cite{wang2022estimation}. 
Hence, 
we may conjecture that the search space is actually smaller and simpler than expected, and thus that a leaner algorithm may produce good results faster. 

In this work, 
we propose an heuristic based on solving a Linear Programming~(LP) model rather than an MILP model at each step of the local search, 
and we generate initial solutions with LP relaxations of the neural surrogate model. 
Confirming our intuition, this strategy is computationally better at scale, such as when neural networks have larger inputs, more neurons, or greater depths.

\section{Notation and Conventions}

In this paper, we consider feedforward networks with fully-connected layers of neurons having ReLU activation. 
Note that convolutional layers can be represented as fully-connected layers with a block-diagonal weight matrix. 
We also abstract that fully-connected layers are often followed by a softmax layer~\cite{bridle1990softmax}, 
since the largest input of softmax matches the largest output of softmax. 

We assume that the neural network has an input $\vx = [\evx_1 ~ \evx_2 ~ \dots ~ \evx_{n_0}]^\top$ from a bounded domain $\sX$ and corresponding output $\vy = [\evy_1 ~ \evy_2 ~ \dots ~ \evy_m]^\top$, and each layer $l \in \sL = \{1,2,\dots,L\}$ has output $\vh^l = [\evh_1^l ~ \evh_2^l \dots \evh_{n_l}^l]^\top$ from neurons indexed by $i \in \sN_l = \{1, 2, \ldots, n_l\}$. 
Let $\mW^l$ be the $n_l \times n_{l-1}$ matrix where each row corresponds to the weights of a neuron of layer $l$, $\mW_i^l$ the $i$-th row of $\mW^l$, and $\vb^l$ the vector of biases associated with the units in layer $l$. With $\vh^0$ for $\vx$ and $\vh^{L}$ for $\vy$,  the output of each unit $i$ in layer $l$ consists of an affine function $\evg_i^l = \mW_{i}^l \vh^{l-1} + \vb_i^l$ followed by the ReLU activation $\evh_i^l = \max\{0, \evg_i^l\}$. 
We denote the neuron \emph{active} when $\evh_i^l = \evg_i^l > 0$ and \emph{inactive} when $\evh_i^l = 0$ and $\evg_i^l < 0$. 
When $\evh_i^l = \evg_i^l = 0$, the state is given by the last nonzero value of $\evg_i^l$ during local search.

In typical neural surrogate models, 
the parameters $\mW^l$ and $\vb^l$ of each layer $l \in \sL$ are constant. 
The decision variables are the inputs of the network ($\vx = \vh^0 \in \sX$) and, in each layer, 
the outputs before and after activation ($\vg^l \in \sR^{n_l}$ and $\vh^l \in \sR_+^{n_l}$ for $l \in \sL$) as well as the activation states ($\vz^l \in \{0,1\}^{n_l}$ for $l \in \sL$).  
By linearly mapping these variables according to the parameters of the network, 
each possible combination of inputs, outputs, and activations become a solution of an MILP formulation. 
For each layer $l \in \sL$ and neuron $i \in \sN_l$, 
the following constraints associate its decision variables $\vh^l$, $\vg^l_i$, $\vh^l_i$, and $\vz^l_i$:
\begin{align}
    \mW_i^l \vh^{l-1} + \vb_i^l = \vg_i^l \label{eq:mip_unit_begin} \\
    (\vz_i^l = 1) \rightarrow \vh_i^l = \vg_i^l \label{eq:first_indicator} \\  
    (\vz_i^l = 0) \rightarrow (\vg_i^l \leq 0 \wedge \vh_i^l = 0) \label{eq:last_indicator} \\
    \vh_i^l \geq 0 \label{eq:lp_unit_end} \\ 
    \vz^l_i \in \{0,1\} \label{eq:mip_unit_end}
\end{align}
The indicator constraints (\ref{eq:first_indicator})--(\ref{eq:last_indicator}) can be modeled with big M constraints~\cite{bonami2015indicator}.

We follow the convention of characterizing each linear region by the set of neurons that they activate~\cite{raghu2017expressive}. For an input $\vx$, 
let $\sS^l(\vx) \subseteq \{ 1, 2, \ldots, n_l\}$ denote the \emph{activation set} of layer $l$.  
Hence, layer $l$ defines an affine transformation of the form $\Omega^{\sS^l(\vx)}(\mW^l \vh^{l-1} + \vb^l)$, where $\Omega^{\sU}$ is a diagonal $v \times v$ matrix in which $\Omega^{\sU}_{ii} = 1$ if $i \in \sU$ and $\Omega^{\sU}_{ii} = 0$ 
otherwise for a subset $\sU \subseteq \sV = \{1, 2, \ldots, v \}$. 
For the linear region containing $\vx = \vx^0$, 
the output of the neural network is the affine transformation $\mT \vt + \vt$ 
for $\mT = \prod_{l=1}^L \Omega^{\sS^l(\vx^0)} \mW^{l}$ 
and $\vt = \sum_{l'=1}^{L} \left( \prod_{l''=l'+1}^{L} \Omega^{\sS^{l''}(\vx^0)} \mW^{l''} \right) \Omega^{\sS^{l'}(\vx^0)} \vb^{l'}$, 
in comparison to which we note that the output $h^{\ell}$ of layer $\ell$ is obtained by replacing $L$ with $\ell$~\cite{huchette2023survey}. 

\section{Walking Along Linear Regions}

Let us 
consider a neural network representing the piecewise linear function $f(\vx)$, a linear objective function $F(\vx) = \vc^\top f(\vx)$ to be maximized, and an implicit set of linear constraints from assuming the input set $\sX$ to be a polytope. 

We can model our problem as an MILP on $\left( \vx, \{ \vg \}_{i=1}^L, \{ \vh \}_{i=0}^L, \{ \vz \}_{i=1}^L, \vy \right)$:
\begin{align}
\max ~ & 
\vc^\top \vy \label{eq:mip_nw_begin} \\ 
\text{s.t.} ~ & 
\text{(\ref{eq:mip_unit_begin})--(\ref{eq:mip_unit_end})} &
\forall l \in \sL, i \in \{1, \ldots, n_l \} \\
& \vx = \vh^0, \vy = \vh^L, \vx \in \sX \label{eq:mip_nw_end}
\end{align}
For an input $\vx = \vx^0$, we can define an LP model by fixing the binary variables as 
$\vz^l_i = 1$ if $i \in \sS^{\ell}(\vx^0)$ and 
$\vz^l_i = 0$ otherwise. 
Let us denote it as $\textbf{LP}(\vx^0)$.
By not fixing $\vx$, $\textbf{LP}(\vx^0)$ finds an input maximizing $F(\vx)$ in a linear region with $\vx^0$.

We propose the local search outlined in Algorithm~\ref{alg:local_search}, 
which is a loop moving from an input $\vx^0$ to the input $\vx^1$ in the same linear region by solving $\text{LP}(\vx^0)$. 
If we find an improvement, we continue moving along the same direction $\vd = \vx^1 - \vx^0$ to the next linear region with a step  $\varepsilon \vd$ updating $\vx^0$. 
We expect that $F(\vx^1 + \varepsilon \vd) > F(\vx^1)$ since $\| \nabla F(\vx^1 + \varepsilon \vd) - \nabla F(\vx^1) \|$ is usually small~\cite{wang2022estimation}. 
Ideally, $\varepsilon$ should be small enough to move only to the next linear region while being large enough 
to be numerically computed as in the relative interior of that next linear region. 
We also adjust the move along each dimension to ensure that $\vx^0 \in \sX$. 

    \begin{algorithm}[t]
    \caption{Local search to walk within and across linear regions.} 
    \label{alg:local_search}
    {\footnotesize
    \begin{algorithmic}[1]
    \Repeat \Comment{Local search consists of an improvement loop}
    \State $\vx^1 \gets$ Optimal solution of $\textbf{LP}(\vx^0)$ \Comment{Finds best solution within linear region}
    \If{$F(\vx^1) > F(\vx^0)$} \Comment{Checks if there was an improvement}
    \State $\vd \gets \vx^1-\vx^0$ \Comment{Computes direction of improvement $\vd$}
    \State $\vx^0 \gets \vx^1 + \varepsilon \vd$
    \Comment{Leaves the linear region along direction $\vd$}
    \For{$i \gets 1, \ldots, n_0$} \Comment{Loops over all input dimensions}
    \If{$x^1 + \varepsilon \vd e^i \notin \sX$} \Comment{Checks if move is outside input space}
    \State $x^0_i \gets x^1_i$ \Comment{Corrects move to be inside input space}
    \EndIf
    \EndFor
    \EndIf
    \Until{$F(\vx^1) \leq F(\vx^0)$} \Comment{Stops when no improvement occurs}
    \State \Return $\vx^0$ \Comment{Returns best solution found}
    \end{algorithmic}
    }
    \end{algorithm}

\begin{figure}
\centering
\includegraphics[width=0.7\textwidth]{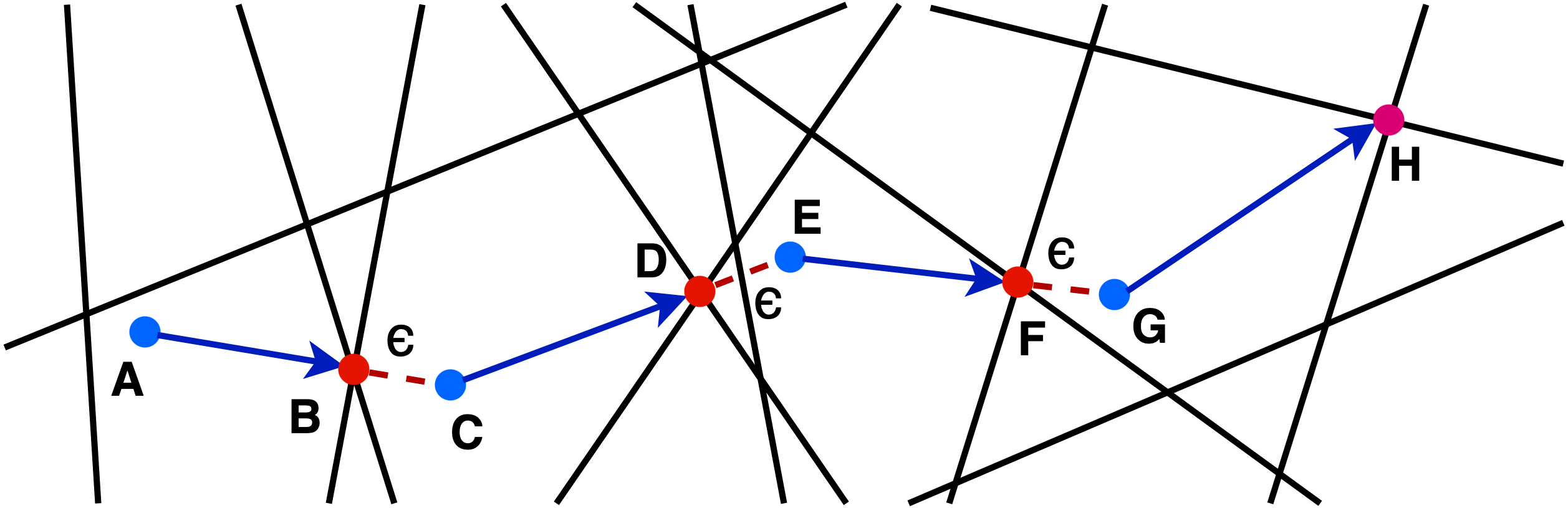}
\caption{From a starting point, our local search algorithm moves in a certain direction indicated by the blue arrow, and then takes a small step into the next linear region before moving again. We stop when the next linear region has no better solution.  
} \label{fig:walk}
\end{figure}

Figure~\ref{fig:walk} illustrates three iterations of improvement with the local search algorithm, 
each characterized by a pair of points $(\vx^0, \vx^1)$ denoting a direction of improvement: $(A,B)$, $(C,D)$, and $(E,F)$. 
Among those, the second iteration shows that a larger step may skip a smaller linear region. 
Conversely, we could mistakenly conclude that no further improvement is possible if a smaller step $\vd$ is numerically computed in such a way that $\vx^1 = H \approxeq \vx^1 + \varepsilon \vd$. 

We embed the local search in a generator of initial solution outlined in Algorithm~\ref{alg:generator}, 
which is based on solving variations of the linear relaxation of (\ref{eq:mip_nw_begin})--(\ref{eq:mip_nw_end}). 
We can compute an input $\tilde{x}$ that is somewhat aligned with maximizing function $F(\vx)$ by solving this relaxation. 
We denote this model as $\text{LR}$:
\begin{align}
\max ~ & 
\vc^\top \vy \label{eq:rl_nw_begin} \\ 
\text{s.t.} ~ & 
\text{(\ref{eq:mip_unit_begin})--(\ref{eq:lp_unit_end})} &
\forall l \in \sL, i \in \{1, \ldots, n_l \} \\
& \vz^l_i \in [0,1] & \forall l \in \sL, i \in \{1, \ldots, n_l \} \\
& \vx = \vh^0, \vy = \vh^L, \vx \in \sX \label{eq:rl_nw_end}
\end{align}
We then impose a random sequence of constraints on activation states, producing a solution $\bar{x}$ in a different linear region after adding each constraint. The probability of fixing a neuron are calibrated to produce the most change to the linear relaxation. We start over from $\tilde{x}$ when no more activations can be fixed. 

\paragraph{Related Work} We denote our approach as Relax-and-Walk (\textbf{RW}) and the local search in~\cite{perakis2022optimizing} as ``Sample-and-MIP'' (\textbf{SM}). SM is based on generating initial solutions by random sampling and then solving a restriction of the MILP (\ref{eq:mip_nw_begin})--(\ref{eq:mip_nw_end}) to identify the best solution among adjacent linear regions. SM may find the best solution locally, 
but it may take much longer to compute in networks with larger dimensions. 
Hence, SM may produce fewer solutions and less improvement.

    \begin{algorithm}[!tb]
    \caption{Generation of initial solutions and injection in local search}
    \label{alg:generator}
    {\footnotesize
    \begin{algorithmic}[1]
    \State $( \widetilde{\vx}, \{ \widetilde{\vg} \}_{i=1}^L, \{ \widetilde{\vh} \}_{i=0}^L, \{ \widetilde{\vz} \}_{i=1}^L, \widetilde{\vy} ) \gets$ Optimal solution of $\text{LR}$
    \State \textbf{DO LOCAL SEARCH FROM $\widetilde{\vx}$} \Comment{First local search; and the only one at $\widetilde{\vx}$}
    \Loop \Comment{Outer loop defines indefinite run until interrupted}
    \State $(\bar{\vx},\bar{\vz}) \gets (\widetilde{\vx}, \widetilde{\vz})$ \Comment{Uses first $\text{LR}$ solution to decide where to go next}
    \For{$\ell \gets 1, \ldots, L$} \Comment{Loops sequentially over layers to fix neurons}
    \State $\sN \gets \{ 1, 2, \ldots, n_\ell \}$  \Comment{Accounts for all neurons to fix in layer $\ell$}
    \While{$\sN \neq \emptyset$} \Comment{Loops to try fixing each neuron once}
    \For{$i \in N$} \Comment{Loops over unfixed neurons}
    \If{$i \in \sS^\ell(\bar{\vx})$} \Comment{Checks if neuron $i$ is active in last solution $\bar{x}$}
    \State $\chi_i \gets 1 - \bar{\vz}_i^\ell$ \Comment{If so, measures distance of relaxed binary to 1}
    \Else
    \State $\chi_i \gets \bar{\vz}_i^\ell$ \Comment{Otherwise, measures distance of relaxed binary to 0}
    \EndIf
    \EndFor \Comment{Produces a shifted probability on $\chi$ values to pick a neuron}
    \State $k \gets$ Element $i \in \sN$ with probability $\chi_i + \delta / \sum_{j \in \sN} (\chi_j + \delta)$
    \State $\sN \gets \sN \setminus \{ k \}$ \Comment{Records attempt to fix neuron $k \in \sN$}
    \If{$k \in \sS^\ell(\bar{x})$} \Comment{Checks if neuron $k$ is active in last solution $\bar{x}$}
    \State Add constraint $\vz_k^\ell = 0$ to $\text{LR}$ \Comment{If so, makes it inactive going forward}
    \Else 
    \State Add constraint $\vz_k^\ell = 1$ to $\text{LR}$ \Comment{Otherwise, makes it active}
    \EndIf
    \If{$\text{LR}$ is feasible} \Comment{Checks if new constraint keeps $\text{LR}$ feasible}
    \State $\left( \bar{\vx}, \{ \bar{\vg} \}_{i=1}^L, \{ \bar{\vh} \}_{i=0}^L, \{ \bar{\vz} \}_{i=1}^L, \bar{\vy} \right) \gets$ Optimal solution of $\text{LR}$
    \State \textbf{DO LOCAL SEARCH FROM $\bar{\vx}$} \Comment{Local search at new solution}
    \Else \Comment{In case not, neuron can only have same activation as before}
    \State Remove constraint on $\vz_k^\ell$; revert $(\bar{x},\bar{z})$ to last feasible solution of $\text{LR}$
    \EndIf
    \EndWhile \Comment{Fixed the entire layer; moves on to the next}
    \EndFor \Comment{Fixed all layers; ready to drop constraints}
    \State Remove all activation constraints from $\text{LR}$
    \EndLoop \Comment{Starts over from $(\widetilde{x}, \widetilde{z})$}
    \end{algorithmic}
    }
    \end{algorithm}



\section{Experiments}

We benchmark our \textbf{RW} method with \textbf{SM}~\cite{perakis2022optimizing} and \textbf{Gurobi} 10.0.1. 
For local search, we use $\varepsilon = 0.01$ since by preliminary tests it was small enough to avoid skipping linear regions. 
We ran the code in \cite{tong2024rw} on 10 cores of a cluster with Intel(R) Xeon(R) Gold 6336Y CPU @ 2.40GHz processors and 16 GB of RAM. 

\subsection{Random ReLU Networks}

Our first experiment replicates and extends the optimization of output value of randomly initialized neural networks in~\cite{perakis2022optimizing} to test scalability and solution quality. 
With a time limit of 1 hour, we use 5 different networks for each choice of input sizes $n_0 \in \{10, 100, 1000\}$ and configurations of the form $L \times n_\ell$ for depth $L \in \{1, 2, 3\}$ and width $n_\ell \in \{100, 500\}$. 
We note that solving to optimality with Gurobi within 1 hour is very unlikely, except for $1 \times 100$ with $n_0=10$.

\begin{figure}[ht]
\includegraphics[width=0.95\textwidth]{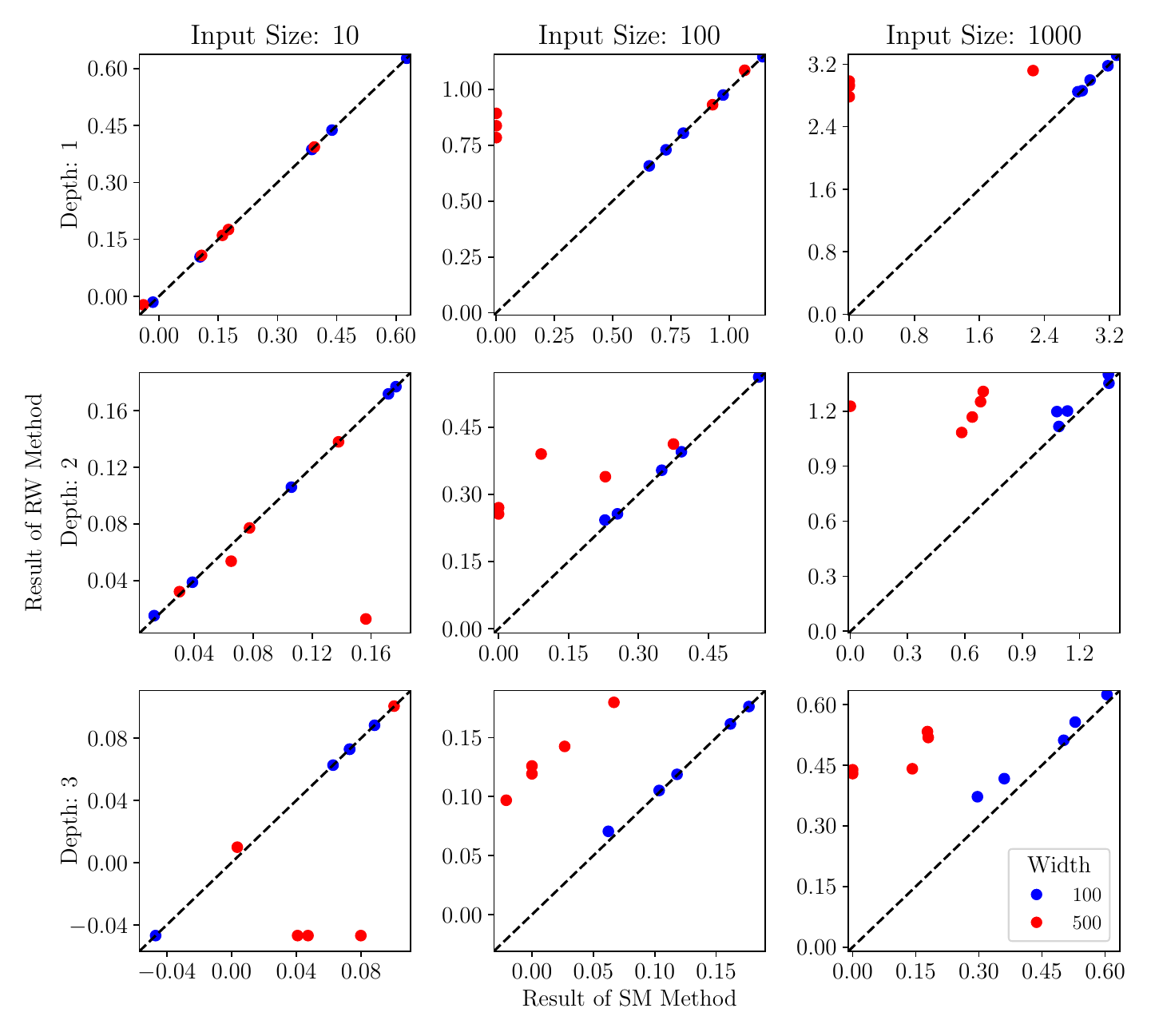}
\vspace{-10pt}
\caption{Comparison of best objective values obtained by RW and SM in random networks.  
The points are favorable to RW above the line $Y=X$; and to SW below it. 
RW is at least 1\% better in 64.4\% of the cases, while SM is in 12.2\%. 
} \label{fig:rw_vs_sm}
\end{figure}

\textbf{RW vs. SM:} Figure~\ref{fig:rw_vs_sm} shows the pair of values obtained for the same random network with RW and SM. RW outperforms SM for $n_0 \in \{ 100, 1000 \}$ and $n_\ell = 500$.
The performance is similar for $n_0 = 10$, except in the four cases where Gurobi fails to solve the linear relaxation. Those have minimum value in the plots. That happens more often when $n_0$ is the smallest while $L$ is larger:  
the model is likely more sensitive to numerical issues as the linear regions get smaller.

\begin{table}[h!]
\caption{Average number of final solutions produced by local search with each method.}
    \label{tab:walk_vs_mip}
\centering
    \begin{tabular}{c@{\hspace{1cm}}cc@{\hspace{0.25cm}}cc@{\hspace{0.5cm}}cc@{\hspace{0.25cm}}cc@{\hspace{0.5cm}}cc@{\hspace{0.25cm}}c}
    \multirow{2}{*}{Model Size} & \multicolumn{3}{c@{}}{$n_0=10$} & & \multicolumn{3}{c@{}}{$n_0=100$} & & 
    \multicolumn{3}{c@{}}{$n_0=1000$} \\
    \cline{2-4}
    \cline{6-8}
    \cline{10-12}
    & RW & & SM & & RW & & SM & & RW & & SM \\
    \midrule
    $1 \times 100$ & 2882.0 & & 2562.2 & & 1038.4 & & 548.0 & & 134.6 & & 100.2 \\
    $1 \times 500$ & 174.8 & & 165.4 & & 80.0 & & 4.4 & & 16.2 & & 1.0 \\
    $2 \times 100$ & 429.2 & & 230.4 & & 331.0 & & 39.8 & & 55.8 & & 11.8 \\
    $2 \times 500$ & 10.0 & & 2.8 & & 25.4 & & 1.0 & & 6.0 & & 1.0 \\
    $3 \times 100$ & 421.8 & & 197.0 & & 204.0 & & 18.2 & & 33.8 & & 4.2 \\
    $3 \times 500$ & 11.5 & & 2.0 & & 15.0 & & 1.0 & & 3.0 & & 1.0 \\
    \end{tabular}
\vspace{-20pt}
\end{table}

Moreover, we observe that \textbf{walking is cheaper than MIPing}: Table~\ref{tab:walk_vs_mip} shows that the walking algorithm RW converges more frequent to a local optimum by the time limit. Conversely, 
the average runtime of MILP restrictions in SM explodes very quickly when the network gets wider, consistent with  Figure~\ref{fig:rw_vs_sm}. 
This can be explained by the number of unfixed MILP variables growing with the network dimensions in the SM approach. 
Moreover, consecutive steps of the SM approach may reevaluate some neighboring linear regions again.

\begin{figure}
\vspace{-10pt}
\includegraphics[width=0.9\textwidth]{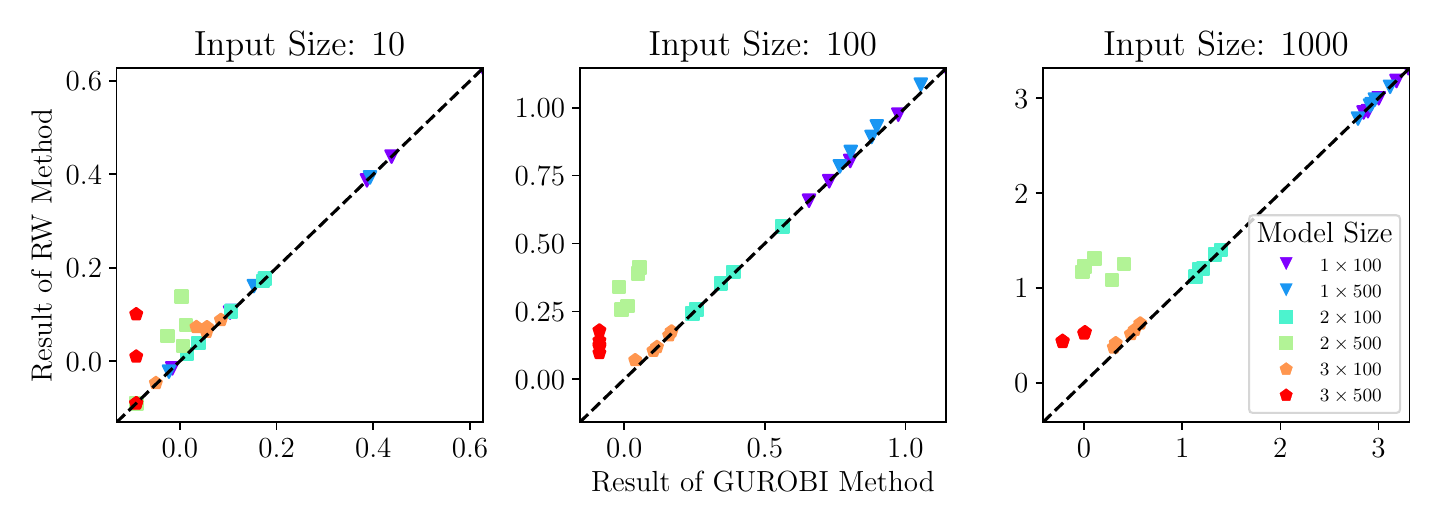}
\vspace{-10pt}
\caption{Comparison of best objective values obtained by RW and Gurobi in random networks. 
RW is at least 1\% better in 55.1\% of the cases, while Gurobi is in 12.4\%.
} \label{fig:rw_vs_gurobi}
\end{figure}

\textbf{RW vs. Gurobi:} Figure~\ref{fig:rw_vs_gurobi} shows a similar comparison between RW and Gurobi, 
but with depths combined for conciseness. 
Gurobi can handle a shallow network ($L=1$) even with the largest input size $n_0=1000$. 
When the network is deeper and the structure of the linear regions more complex~\cite{serra2018bounding}, 
directly solving an MILP is slower. 
When both width and depth are large, Gurobi cannot find a feasible solution---see points next to $Y$-axis. 
The same four cases with unbounded relaxation are also difficult for Gurobi---see points near the left bottom.

\vspace{-25pt}

\subsection{Optimal Adversary Experiment}

~
\begin{wrapfigure}{r}{0.42\textwidth}
  \centering
  \vspace{-60pt}
    \includegraphics[width=0.4\textwidth]{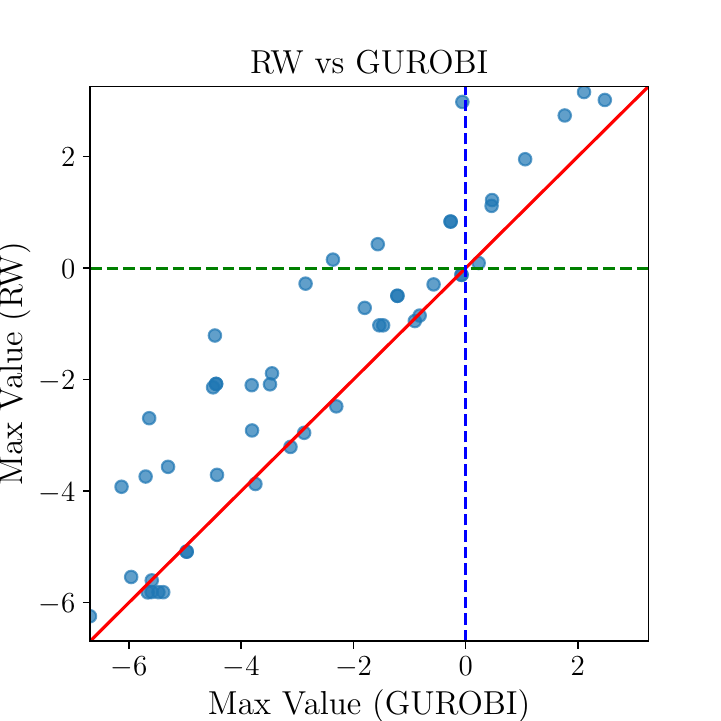}
  \caption{Comparison of best objective values obtained by RW and Gurobi in optimal adversary models. 
RW is at least 1\% better in 68\% of the cases, while Gurobi is in 30\%. }
\label{fig:mnist}
\vspace{-15pt}
\end{wrapfigure}
Given an input $\vx = \hat{\vx}$, 
the output neuron for the predicted label $c$, 
and the output neuron for another likely label $w$, 
the optimal adversary problem aims to maximize $\vy_w - \vy_c$ for an input $\vx$ sufficiently near $\hat{\vx}$. 
We solved this problem for $|\vx - \hat{\vx}|_1 < \Delta$ as in~\cite{tsay2021partition} by using a setup derived from the Gurobi Machine Learning repository~\cite{gurobi2023ml}. 
Figure~\ref{fig:mnist} shows the result from testing $50$ images from the MNIST dataset~\cite{lecun1998mnist} with $\Delta = 5$ for $1$ hour, 
all of which on a $2 \times 500$ classifier with test accuracy $97.04\%$.  
In 10\% of the cases, RW found an adversarial input (positive solution) and Gurobi did not. 
When RW does better, it does so by a wider margin.

\section{Conclusion}

We introduced a local search algorithm for optimizing over trained neural networks. 
We designed our algorithm to leverage model structure based on what is known about linear regions in deep learning.  
Moreover, our algorithm scales more easily because it only solves LP models at every step. 
Last, but certainly not least, the solutions are usually better in comparison with other methods.

\paragraph{Acknowledgement} The authors were supported by 
the National Science Foundation (NSF) grant IIS 2104583, 
including Junyang Cai while at Bucknell.

%
%
%
\bibliographystyle{splncs04}
\bibliography{references}
\end{document}